\documentclass[12pt]{amsart}

\makeatletter
\@addtoreset{equation}{section}
\makeatother

\usepackage{amsmath}
\usepackage{latexsym}
\usepackage{amssymb}
\usepackage{amsfonts}

\usepackage{graphics}

\renewcommand{\proof}{\par\noindent{\it Proof.\ \ }}
\def\qed{\ifmmode\square\else\nolinebreak\hfill
$\square$\fi\par\vskip12pt}

\def\l{\langle} \def\r{\rangle}

 \def\ZZ{\mathbb Z}
\def\BB{{\mathcal B}}

\def\mod{{\sf mod~}}

\def\Aut{{\sf Aut}}  
\def\Cos{{\sf Cos}}  \def\Cay{{\sf Cay}}

\def\K{{\sf K}} 
\bf

\def\Cay{{\rm Cay}}
\def\D{{\rm D}}  
 \def\S{{\rm S}} \def\G{{\rm G}}
\def\J{{\rm J}} \def\M{{\rm M}}
\def\soc{{\rm soc}} 
\def\C{{\bf C}}\def\N{{\bf N}} 

\def\Ga{{\it\Gamma}}  

\def\a{\alpha} \def\b{\beta} \def\d{\delta}
\def\g{\gamma} \def\s{\sigma}

\def\POmega{{\rm P\Omega}}
\def\Sp{{\rm Sp}} 

\def\AGammaL{{\rm A\Gamma L}}
\def\PGammaL{{\rm P\Gamma L}} 
\def\A{{\rm A}}

\def\PSL{{\rm PSL}}\def\PGL{{\rm PGL}}

\def\AGL{{\rm AGL}}\def\ASL{{\rm ASL}}
 \def\PSU{{\rm PSU}}

  \def\D{{\rm D}} \def\G{{\rm G}}

\newtheorem{theorem}{Theorem}[section]%
\newtheorem{lemma}[theorem]{Lemma}%
\newtheorem{corollary}[theorem]{Corollary}%
\newtheorem{example}[theorem]{Example}%
\newtheorem{remark}[theorem]{Remark}%

\begin{document}

\title[Edge-transitive graphs]{On the automorphism groups \\ of
graphs with twice  prime valency}
\thanks{This work was 
partially supported by the National Natural Science
Foundation of China (11731002) and the Fundamental Research Funds for the Central Universities.}
\thanks{Corresponding author: Zai Ping Lu (lu@nankai.edu.cn)}

\author[Liao]{Hong Ci Liao}
\address{H. C. Liao\\ Center for Combinatorics\\
LPMC-TJKLC, Nankai University\\
Tianjin 300071\\
P. R. China} \email{827436562@qq.com}

\author[Li]{Jing Jian Li}
\address{Jing Jian Li\\ School of Mathematics and Information Sciences\\ Guangxi
University\\ Nanning 530004, P. R. China.}
\address{Colleges and Universities
Key Laboratory of Mathematics and Its Applications.}
\email{lijjhx@163.com}

\author[Lu]{Zai Ping Lu}
\address{Z. P. Lu\\ Center for Combinatorics\\
LPMC-TJKLC, Nankai University\\
Tianjin 300071\\
P. R. China} \email{lu@nankai.edu.cn}
\maketitle
\date\today

\begin{abstract}
A graph is edge-transitive  if its automorphism group acts transitively on the edge set. In this  paper, we investigate the automorphism groups of edge-transitive graphs of odd order and twice prime valency.
Let $\Ga$ be a connected  graph of odd order and twice prime valency, and let $G$ be a subgroup of the automorphism group of $\Ga$.
In the case where $G$ acts transitively on the edges and quasiprimitively on the vertices of $\Ga$, we prove that either $G$ is almost simple or $G$ is a primitive group of  affine type. If further $G$ is an almost simple primitive group then, with two exceptions, the socle of $G$ acts transitively on the edges of $\Ga$.

\vskip 10pt

\end{abstract}

\vskip 50pt
\section{introduction}

In this paper, all graphs are assumed to be finite and simple. In particular, a graph  is a pair $\Ga=(V,E)$
of a nonempty set $V$ and a set $E$ of $2$-subsets of $V$, which are called the vertex set and edge set of $\Ga$, respectively.
Each edge $\{\a,\b\}\in E$ gives two ordered pairs $(\a,\b)$ and $(\b,\a)$, every of them is called an arc of $\Ga$.
A triple $(\a,\b,\g)$ of vertices is a $2$-arc if $\a\ne \g$ and both $(\a,\b)$ and $(\b,\g)$ are arcs. 

\vskip 5pt

Assume that $\Ga=(V,E)$ is a graph. An automorphism $g$ of $\Ga$ is a permutation (i.e., a bijection) on $V$ such that $\{\a^g,\b^g\}\in E$ for all $\{\a,\b\}\in E$.
Denote by $\Aut\Ga$ the set of all automorphisms of $\Ga$. Then $\Aut\Ga$ is a (finite) group under the product of permutations, which acts naturally on the edge set, arc set and $2$-arc set of $\Ga$ by
\[\{\a,\b\}^g=\{\a^g,\b^g\},\, (\a,\b)^g=(\a^g,\b^g),\, (\a,\b,\g)^g=(\a^g,\b^g,\g^g),\]
respectively. For a subgroup $G\le \Aut\Ga$, the graph $\Ga$ is said to be
$G$-vertex-transitive, $G$-edge-transitive,  $G$-arc-transitive and  $(G,2)$-arc-transitive if $G$ acts transitively on the vertex set, edge set, arc set and $2$-arc set of $\Ga$, respectively.

\vskip 5pt

Let $\Ga=(V,E)$ be a connected graph, and $G\le \Aut\Ga$. Suppose that $G$ acts quasiprimitively on $V$, that is, every minimal normal subgroup of $G$ has a unique orbit on $V$.
Following the subdivision  in \cite{Prag-quasi},
the group $G$ is one of eight types of quasiprimitive groups.
Conversely, using the `coset graph' construction (see \cite{Sabiddusi}), one can obtain  graphs from each type of quasiprimitive groups. Nevertheless, it is believed that the group $G$ is  quite restrictive if further the graph $\Ga$ is assumed  with certain symmetric properties, or some restrictions on the order or valency.
For example, if $\Ga$ is $(G,2)$-arc-transitive then
Praeger \cite{Prag-o'Nan} proved that only four of those $8$ types occur for $G$. If $\Ga$ has odd order and  $\Ga$ is $(G,2)$-arc-transitive, then $G$ is an almost simple group by \cite{Li-odd}.
In   this paper, replacing the `$2$-arc-transitivity'  by `edge-transitivity' and a restriction on the valency of graph, we  investigate the pair of $\Ga$ and $G$.

\vskip 5pt

Let $\Ga=(V,E)$ be a connected graph of twice valency, and $G\le \Aut\Ga$. In \cite{Prag-Xu}, Praeger and Xu give a nice characterization for the graph $\Ga$ while it is  $G$-arc-transitive and $G$ contains an irregular abelian (and so, intransitive) normal subgroup. In Section 3 of this paper,
we focus on the case where $\Ga$ is  $G$-edge-transitive and $G$ contains a
transitive normal subgroup. In particular, letting $\soc(G)$ be the socle of $G$, the following result is proved.

\begin{theorem}\label{myth}
Let $\Ga=(V,E)$ be a connected graph of odd order and valency $2r$ for some prime $r$, and let $G\le \Aut\Ga$. Assume that $\Ga$ is  $G$-edge-transitive but not $(G,2)$-arc-transitive. If $G$ is quasiprimitive on $V$ then either
$G$ is almost simple, or  $\soc(G)=\ZZ_p^k$ for some odd prime $p$  and integer $1\le k\le r$.
\end{theorem}

\vskip 5pt

\begin{remark}\label{rem}
Li and the last two authors of this paper have been working on the project
classifying all graphs of odd order which admit an almost simple group $X$ acting $2$-arc-transitively. By their classification, such a graph $\Ga=(V,E)$ is not of twice prime valency unless one of the following holds: (1) $\Ga$ is a complete graph, (2) $\Ga$ is the odd graph of valency $4$, (3) $\soc(X)=\PSL(2,q)$ and $\Ga$ is of valency $4$ or $6$. In view this, if $X$ has a primitive subgroup acting transitively on $E$ then we may easily prove that $\Ga$ is a complete graph, see also the following result.
\end{remark}

\vskip 5pt

 A graph $\Ga$ is called $2$-arc-transitive if it is $(\Aut\Ga,2)$-arc-transitive. The following result is a consequence of Theorem \ref{myth}, which is proved in Section 4.

\begin{theorem}\label{myth-2}
Let $\Ga=(V,E)$ be a connected graph of odd order and twice prime valency, and let $G\le \Aut\Ga$. Assume that $\Ga$ is $G$-edge-transitive but not $(G,2)$-arc-transitive. If $G$ is primitive on $V$ then either $\Ga$ is a complete graph, or  $\soc(\Aut\Ga)=\soc(G)$; in particular,
$\Ga$ is $2$-arc-transitive
if and only if $\Ga$ is a complete graph.
\end{theorem}

For a pair $(G,\Ga)$ in Theorem \ref{myth} or \ref{myth-2},  the action of   $\soc(G)$ on the edge set of $\Ga$ is considered in Section 5. We present several examples and prove the following result.

\begin{theorem}\label{myth-3}
Let $\Ga=(V,E)$ be a connected graph of odd order and twice prime valency, and let $G\le \Aut\Ga$. Assume that $\Ga$ is $G$-edge-transitive but not $(G,2)$-arc-transitive. If $G$ is almost simple and primitive on $V$ then either $\Ga$ is $\soc(G)$-edge-transitive, or $\Ga$ is of valency $4$
and isomorphic one of the graphs in Example \ref{exam-2}.
\end{theorem}

\vskip 30pt

\section{Preliminaries}
For convenience,  a graph $\Ga=(V,E)$ is sometimes viewed as the digraph on $V$ with directed edges the arcs of $\Ga$, and a subset $\Delta$ of the arc set of $\Ga$ is sometimes viewed as a digraph on $V$.
Then, with such a convenience, a digraph  $\Delta$ is a graph if and only if $\Delta=\Delta^*:=\{(\a,\b)\mid (\b,\a)\in \Delta\}$.
In this section, we make the following assumption:
\begin{enumerate}
\item[] {\em $\Ga=(V,E)$ is a connected  regular graph,  $G\le \Aut\Ga$, $N\ne 1$ is a normal subgroup of $G$, and $\Delta$ is the union of some $G$-orbits on the arc set of $\Ga$.}
\end{enumerate}

\vskip 5pt

Let $\a\in V$. Set $\Delta(\a)=\{\b\in V\mid (\a,\b)\in \Delta \}$ and $N_\a=\{g\in N\mid \a^g=\a\}$, called respectively the (out-)neighborhood of $\a$ in (the digraph) $\Delta$ and the stabilizer of $\a$ in $N$. Then $N_\a$ induces a  permutation group $N_\a^{\Delta(\a)}$.
Denote by $N_\a^{[1]}$ the kernel of $N_\a$ acting on $\Delta(\a)$.
Then
\begin{equation}\label{eq-1}
N_\a^{\Delta(\a)}\cong N_\a/N_\a^{[1]},\, N_{\a^g}=N_\a^g,\, N_{\a^g}^{[1]}=(N_\a^{[1]})^g,\, \forall g\in G.
\end{equation}
For   $U\subseteq  V$, set $N_{(U)}=\cap_{\b\in U}N_\b$ and $N_U=\{x\in N\mid U^x=U\}$, called respectively
the point-wise stabilizer and set-wise stabilizer of $U$ in $N$. Then
\begin{equation}\label{eq-2}
N_{(U^g)}=N_{(U)}^g, \, N_{U^g}=N_U^g,\, \forall g\in G.
\end{equation}
If $U=\{\a_1,\ldots,\a_k\}$ then we write $N_{(U)}=N_{\a_1\a_2\ldots\a_k}$.

\vskip 5pt

For a finite group $X$, denote by $\pi(X)$   the set of prime divisors of $|X|$.
The following lemma says that $\pi(N_\a)=\pi(N_\a^{\Delta(\a)})$ if $G$ is transitive on $V$, see also \cite[Lemma 1.1]{CLP2000} for the case where $\Delta$ is the arc set of $\Ga$.

\begin{lemma}\label{normal-subg-1}
Assume that $\Delta$ is connected and $G$ is transitive on $V$. Then
$\pi(N_\a)=\pi(N_\a^{\Delta(\a)})$ and $|\Delta(\a)|\ge \max\pi(N_{\a})$,
 where $\b\in \Delta(\a)$. If further $\Delta=\Delta^*$ then $|\Delta(\a)|> \max\pi(N_{\a\b})$.
\end{lemma}
\proof
Let $p$ be an arbitrary prime. Suppose that $p$ is a not divisor of $|N_\a^{\Delta(\a)}|$. Then $P\le N_\a^{[1]}$, where   $P$ is an arbitrary Sylow $p$-subgroup of $N_\a$. In particular,
$P\le N_\b$ for $\b\in \Delta(\a)$.
Take $g\in G$ with $\a^g=\b$. We have $\Ga(\b)=\Ga(\a)^g$,
$N_\b=N_\a^g$ and $N_\b^{[1]}=(N_\a^{[1]})^g$. It follows that  $P$ is also a Sylow $p$-subgroup of $N_\b$, and $P\le N_\b^{[1]}$.

 Let $\g\in V\setminus\{\a,\b\}$. By \cite[Lemma 2]{Neumann}, there are $\b=\b_0, \b_1, \ldots\b_t=\g$ such that $(\b_{i-1},\b_{i})\in \Delta$ for $1\le i\le t$. By the argument in the previous paragraph and induction on $i$, we conclude that $P\le N_\g^{[1]}$. It follows that  $P\le N_{(V)}=1$, that is, $p\not\in \pi(N_\a)$. Then we have  $\pi(N_\a)=\pi(N_\a^{\Delta(\a)})$.

Let $q$ be a prime with $q>|\Delta(\a)|$, and let $Q$ be a Sylow $q$-subgroup of $N_{\a}$. Then $Q$
 $Q$ acts trivially on $\Delta(\a)$.
 Since $|\Delta(\b)|=|\Delta(\a)|<q$, we conclude that $Q$ acts trivially on $\Delta(\b)$.
For $\g\in V\setminus\{\a,\b\}$,
choose $\b=\b_0, \b_1, \ldots\b_t=\g$ such that $(\b_{i-1},\b_{i})\in \Delta$ for $1\le i\le t$. Noting that $|\Delta(\b_i)|=|\Delta(\a)|$ for each $i$, by induction on $i$, we have $Q\le N_\g$. Thus  $Q\le N_{(V)}=1$, and so $q\not\in \pi(N_{\a})$.
Then $|\Delta(\a)|\ge \max\pi(N_{\a})$.

Finally, by \cite[Lemma 1.1]{CLP2000}, if $\Delta=\Delta^*$ then $|\Delta(\a)|>\max\pi(G_{\a\b})\ge \max\pi(N_{\a\b})$. Thus  the lemma follows.
\qed

\begin{lemma}\label{normal-subg-2}
Assume that $G$ is transitive on   $V$ and $\Delta$ is a $G$-orbit. Suppose that $\Delta$ is connected.
Then
$\pi(N_{\a\b}^{\Delta(\b)})=\pi(N_{\a\b})$, where $\b\in \Delta(\a)$.
\end{lemma}
\proof
Let $p$ be a prime with $p\not\in \pi(N_{\a\b}^{\Delta(\b)})$. Then $P\le N_\b^{[1]}$, where   $P$ is an arbitrary Sylow $p$-subgroup of $N_{\a\b}$. In particular, $P\le N_{\b\d}$ for $\d\in \Delta(\b)$. Choose $g\in G$ with $(\a,\b)^g=(\b,\d)$. Then $N_{\b\d}=N_{\a\b}^g$, $\Delta(\d)=\Delta(\b)^g$ and $N_{\d}^{[1]}= (N_{\b}^{[1]})^g$, which yields that $P$ is a  Sylow $p$-subgroup of  $N_{\b\d}$, and $P\le N_{\d}^{[1]}$. For $\g\in V\setminus\{\a,\b\}$,
choose $\b=\b_0, \b_1, \ldots\b_t=\g$ such that $(\b_{i-1},\b_{i})\in \Delta$ for $1\le i\le t$. By induction on $i$, we may prove $P\le N_{\g}^{[1]}$.
It implies that $P=1$, so $p\not \in \pi(N_{\a\b})$. Then the lemma follows.
\qed

Suppose that $N$ is intransitive on $V$. Let $\BB$ be the set of $N$-orbits on $V$.
Define a digraph $\Delta_N$ on $\BB$ such that
$(B,C)$ is an arc if and only of $B\ne C$ and $(\d,\g)\in \Delta$ for some
$\d\in B$ and some $\g\in C$.
For the case where $\Delta$ is the arc set of $\Ga$, we also write $\Delta_N$ as $\Ga_N$. Let $B\in \BB$ and $\a\in B$. Then
\begin{equation}\label{eq-3}
\Delta(\a)=(B\cap \Delta(\a))\cup (\cup_{C\in \Delta_N(B)}(C\cap \Delta(\a)),\,|\Delta(\a)|\ge |B\cap \Delta(\a)|+|\Delta_N(B)|.
\end{equation}
The following lemma is easily shown.

\begin{lemma}\label{quotient}
Assume that  $\Delta$ is connected and a $G$-orbit. If $N$ is intransitive on $V$, then $|\Delta(\a)|=|\Delta_N(B)||C\cap \Delta(\a)|$ for each given $C\in \Delta_N(B)$; in this case, $G$ induced a group acting transitively on the arcs of $\Delta_N$.
\end{lemma}

\begin{lemma}\label{N-semiregular}
Assume that $G$ is transitive on $V$,  $N$ intransitive on $V$ and $\Delta$ is connected. If $|\Delta(\a)|=|\Delta_N(B)|$ for some $N$-orbit $B$ and $\a\in B$,
then $N$ is semiregular on $V$, that is, $N_\a=1$.
\end{lemma}
\proof
Assume that $|\Delta(\a)|=|\Delta_N(B)|$. Then $\Delta(\a)\cap B=\emptyset$, and
$|\Delta(\a)\cap C|=1$ for every $C\in \Delta_N(B)$. It follows that $N_\a$ fixes $\Delta(\a)$ point-wise, so $N_\a^{\Delta(\a)}=1$. Then $N_\a=1$ by Lemma \ref{normal-subg-1}.
\qed

\vskip 5pt

Suppose that $\Ga=(V,E)$ is  connected and $G$-arc-transitive. For an edge $\{\a,\b\}\in E$, take $x\in G$ with $(\a,\b)^x=(\b,\a)$. Then $x\in \N_G(G_{\a\b})$ and, since $\Ga$ is connected, $G=\l x,G_\a\r$. Further, such an $x$ may be chosen as a $2$-element in $\N_G(G_{\a\b})$. Thus following lemma holds.

\begin{lemma}\label{coset}
If $\Ga=(V,E)$ is connected and $G$-arc-transitive. Then, for  $\{\a,\b\}\in E$, there is a $2$-element $x$ in $\N_G(G_{\a\b})$ with $x^2\in G_{\a\b}$ and $G=\l x,G_\a\r$.
\end{lemma}

\vskip 5pt

Suppose that $G$ has regular subgroup   $R$. Then, fixing a vertex $\a$, we have a bijection
\[\theta: V\rightarrow R,\, \a^x\mapsto x.\]
Set $S=\{x\in R\mid \a^x\in \Delta(\a)\}$. Then $\theta$ gives an isomorphism from $\Delta$ to the Cayley digraph $\Cay(R,S)$, which is defined on $R$ such
$(x,y)$ is an arc if and only if $yx^{-1}\in S$. Under the isomorphism $\theta$,
the group $\N_{\Aut\Delta}(R)$ corresponds to the normalizer   $\hat{R}{:}\Aut(R,S)$ of $\hat{R}$ in $\Aut\Cay(R,S)$, where
$\hat{R}$ is the image of the action of $R$ by right multiplication on $R$, and
\[\Aut(R,S)=\{\s\in \Aut(R)\mid S^\s=S\}.\]
Clearly, $S$ does not contains the identify element of $R$, and $\Delta$ (or $\Cay(R,S)$) is a graph if and only if $S=S^{-1}:=\{x^{-1}\mid x\in S\}$.
The following lemma collects several well-known facts on Cayley digraphs.

\begin{lemma}\label{cay}
\begin{enumerate}
\item[(1)]  $\Delta$ is isomorphic a Cayley digraph
$\Cay(R,S)$ if and only if $\Aut\Delta$ has a regular subgroup isomorphic to $R$.
\item[(2)] A Cayley digraph $\Cay(R,S)$ is connected if and only if
$\l S\r=R$, and $\Cay(R,S)$ is a graph if and only if $S=S^{-1}$.
\item [(3)] $\N_{\Aut\Cay(R,S)}(\hat{R})=\hat{R}{:}\Aut(R,S)$, and if
$\l S\r=R$ then $\Aut(R,S)$ acts faithfully on $S$.
\end{enumerate}
\end{lemma}


\vskip 30pt

\section{Proof of Theorem \ref{myth}}
In this section, we  assume that
$\Ga=(V,E)$ is a connected    graph of valency $2r$ for some prime $r$, and  $G\le \Aut\Ga$ such that $G$ is transitive on both $V$ and $E$ 
but not transitive on the $2$-arc set of $\Ga$.

\begin{lemma}\label{p-divisor-arc}
Suppose that $\Ga$ is $G$-arc-transitive.
Let $N\ne 1$ be a normal subgroup of $G$. Then
either $N$ is transitive on $V$ or $G$ induces an arc-transitive group on $\Ga_N$. Moreover, one of the following holds.
\begin{enumerate}
\item[(1)] $N$ is semiregular on $V$.
\item[(2)] $N_\a$ is a nontrivial $2$-group, and either $\Ga_N$ has valency $r$ or $N$ has at most $2$ orbits on $V$.
\item[(3)] $r=\max\pi(N_\a)$ is odd, and
\begin{enumerate}
\item[(i)] $\Ga_N$ is a cycle and $G$ induces the full automorphism group of this cycle; or
\item[(ii)] $\Ga$ is bipartite with two parts $N$-orbits; or
\item[(iii)]$N$ is transitive on both $V$ and $E$; or
\item[(iv)] $E$ has a partition $E=E_1\cup E_2$ such that
$(E_1,E_2)^g=(E_2,E_1)$ for some $g\in G\setminus N$, and $(V,E_i)$ are $N$-arc-transitive graphs of valency $r$.
\end{enumerate}
\end{enumerate}
\end{lemma}
\proof
The first part of this lemma follows from Lemma \ref{quotient}. We next show that one of (1)-(3) occurs.
By Lemma \ref{normal-subg-1}, $\pi(N_\a)=\pi(N_\a^{\Ga(\a)})$.
 If $N_\a^{\Ga(\a)}=1$ then $N_\a=1$, and hence $N$ is semiregular, so (1) follows.
In the following, we let $N_\a^{\Ga(\a)}\ne 1$.

  Since $N_\a^{\Ga(\a)}$ is normal in $G_\a^{\Ga(\a)}$, all $N_\a^{\Ga(\a)}$-orbits on $\Ga(\a)$ have the same length
say $l$. Then $l=2,\,r$ or $2r$.
If $l=2$ then $N_\a^{\Ga(\a)}$ is a $2$-group, and so does $N_\a$;
if $r=2$ and $l=2r=4$ then, since $\Ga$ is not $(G,2)$-arc-transitive, $N_\a^{\Ga(\a)}\le G_\a^{\Ga(\a)}\le \D_8$, so $N_\a$ is a $2$-group. These two case give (2) by Lemma \ref{quotient}.
Thus, we next let $l>2$ and $r$ be odd.

Assume first that $G_\a^{\Ga(\a)}$ is primitive. Then $N_\a^{\Ga(\a)}$
is a transitive normal subgroup of  $G_\a^{\Ga(\a)}$. It follows that
$N$ is transitive on $E$ and has at most two orbits on $V$.
Since $G_\a^{\Ga(\a)}$ is not $2$-transitive, by \cite[Theorem 1.51]{Gorenstein}, $r=5$ and $G_\a^{\Ga(\a)}\cong \A_5$ or $\S_5$.
Then  $r=5=\max\pi(N_\a)$, and so (ii) or (iii) holds.

Now assume   that $G_\a^{\Ga(\a)}$ is imprimitive. Then $G_\a^{\Ga(\a)}\le \S_2\wr\S_r$ or $\S_r\wr\S_2$, yielding
$r=\max\pi(G_\a)$. Recalling that $l=2r$ or $r$, we have $r=\max\pi(N_\a)$.
If $l=2r$ then $N_\a^{\Ga(\a)}$ is transitive on $\Ga(\a)$, which yields that
$N$ is transitive on $E$, and  so (ii) or (iii) holds.

Let $l=r$. If $N$ has at three orbits on $V$ then (i) follows from Lemma \ref{quotient}. If $N$ has exactly two orbits on $V$ then $\Ga$ is bipartite,
and (ii) occurs.
Suppose that $N$ is transitive on $V$.
Then $N$ has exactly  two orbits on  the arc set of $\Ga$, say $\Delta_1$ and $\Delta_2$, and either $\Delta_2=\Delta_1^*\ne \Delta_1$ or $\Delta_i=\Delta_i^*$ for each $i$.
These two cases yield parts (iii) and (iv), respectively.
\qed

\begin{corollary}\label{odd-order-p-divisor-arc}
Suppose that $\Ga$ is $G$-arc-transitive.
Let $N$ be a transitive normal subgroup of $G$. 
Then one of the following holds.
\begin{enumerate}
\item[(1)] 
 Either  $N_\a$ is a $2$-group, or
$N$ is transitive on $E$;
\item[(2)]
 $E$ has a partition $E=E_1\cup E_2$ such that
$(E_1,E_2)^g=(E_2,E_1)$ for some $g\in G\setminus N$, and $(V,E_i)$ are $N$-arc-transitive graphs of valency $r$.
\end{enumerate}
\end{corollary}

\vskip 5pt

\begin{lemma}\label{p-divisor-half}
Suppose that $G$ is intransitive on the arc set of $\Ga$.
Let $N\ne 1$ be a normal subgroup of $G$. Then one of the following holds.
\begin{enumerate}
\item[(1)] $N$ is semiregular on $V$.
\item[(2)] $r=\max \pi(N_\a)$ for $\a\in V$, and either
\begin{enumerate}
\item[(i)] $N$ is transitive on both $V$ and $E$; or
\item[(ii)] $G/K\cong \ZZ_l$, where $l\ge 2$ is the number of $N$-orbits on
$V$, and $K$ is the kernel of $G$ acting on the set of $N$-orbits.
\end{enumerate}
\end{enumerate}
\end{lemma}
\proof
Let $\Delta$ be a $G$-orbit on the arc set of $\Ga$. Then $|\Delta(\a)|=r$, and so $r=\max \pi(G_\a^{\Delta(\a)})=\max \pi(G_\a)$ by Lemma \ref{normal-subg-1}.
Since $N_\a^{\Delta(\a)}$ is normal subgroup of $G_\a^{\Delta(\a)}$ and $G_\a^{\Delta(\a)}$ is transitive, all $N_\a^{\Delta(\a)}$-orbits (on $\Delta(\a)$) have the same length. Then either $N_\a^{\Delta(\a)}=1$, or
$N_\a^{\Delta(\a)}$ is transitive on $\Delta(\a)$. By Lemma \ref{normal-subg-1}, the former case yields
$N_\a=1$, and so $N$ is semiregular.

Assume that $N_\a^{\Delta(\a)}$ is transitive on $\Delta(\a)$. Then
$r=\max \pi(N_\a)$. If $N$ is transitive on $V$ then
$N$ is transitive on both $V$ and $E$, and so part (i) occurs.
Suppose that $N$ is intransitive on $V$.
By Lemma \ref{quotient}, $\Delta_N$ is a directed cycle, and $G$ induces
an arc-transitive group on this directed cycle. (Note that $\Delta_N$
has length $2$ if $N$ has two orbits on $V$.) Then part (ii) follows.
\qed

\begin{corollary}\label{odd-order-p-divisor-half}
Suppose that $G$ is intransitive on the arc set of $\Ga$.
Let $N$ be a transitive normal subgroup of $G$. Then  either  $N$ is regular on $V$, or $r=\max \pi(N_\a)$   and $N$ is transitive on both $V$ and  $E$.
\end{corollary}

\vskip 5pt

Recall  that every minimal normal subgroup of a quasiprimtive group
is transitive. Then the following  result finishes
 the proof of Theorem \ref{myth}.

\begin{theorem}\label{odd-order-prim-2r}
Let $N$ be a minimal normal subgroup of $G$. Suppose that  $|V|$ is odd and
$N$ is transitive on $V$. Then $N$ is the unique transitive minimal normal subgroup of $G$, and one of the following holds.
\begin{enumerate}
\item[(1)] $N\cong \ZZ_p^k$ for some odd prime $p$ and positive integer $k\le r$, and $G_\a$ is isomorphic to a subgroup of $\S_2\wr\S_r$; if further $G$ is intransitive on the arc set of $\Ga$, then $k<r$, and $G_\a$ is isomorphic to a subgroup of $\S_r$.
\item[(2)]
$N$ is nonabelian simple.
\end{enumerate}
\end{theorem}
\proof
Write $N=T_1\times\cdots\times T_k$, where $T_i$ are isomorphic simple groups.
Let $M$ be a transitive minimal normal subgroup   of $G$. Suppose that  $M\cap N=1$. Then $M$ centralizes $N$. By \cite[Theorem 4.2A]{Dixon}, both $M$ and $N$ are regular on $V$,  and so $N$ is abelain as $|V|$ is odd. Then we have $M=N$, yielding $M=N=1$, a contradiction. Thus $M\cap N\ne 1$, and so $M=N$. Thus  $N$ is the unique transitive minimal normal subgroup of $G$.

\vskip 5pt

{\bf Case 1}. Assume   that $N$ is soluble. Then $N\cong\ZZ_p^k$ for a prime $p$ and some integer $k\ge 1$. In particular, $N$ is abelian, and so $N$ is regular on $V$, and $p$ is odd.
By Lemma \ref{cay}, writing $\Ga=\Cay(N,S)$, we have that
$G_\a\le \Aut(N,S)$  which acts faithfully on $S$, where $\a$ is the vertex of $\Ga$ corresponding to the identity of $N$.
Since $p$ is odd and $\Ga$ is connected, $N$ is generated by
a half number of the elements in $S$, and thus  $k\le r$.

Note that $\{\{x,x^{-1}\}\mid x\in S\}$ is a $G_\a$-invariant partition of $S$.
It follows that $G_\a$ is isomorphic to a subgroup of $\S_2\wr\S_r$.
Suppose that $\Ga$ is not $G$-arc-transitive. Then $G_\a$ has two orbits on
$S$, say $\{x_1,\ldots,x_r\}$ and $\{x_1^{-1},\ldots,x_r^{-1}\}$, see also \cite[Lemma 3.2]{LLiuL}. This yields that $G_\a$ is isomorphic to a subgroup of $\S_r$.  Consider the element $x=\prod_{i=1}^rx_r$ of $N$. It is easily shown that $G_\a$ centralizes $x$. Then $\l x\r$ is a normal subgroup of $G$. Since $x\in N$, by the minimum of $N$, we have either $k=1$ or $x=1$. If $x=1$ then $N=\l x_1,\ldots, x_{r-1}\r$, and thus $k<r$. Thus part (1) follows.

\vskip 5pt

{\bf Case 2}. Assume that $T_i$ are isomorphic nonabelian simple groups.
Since $|V|=|N:N_\a|$ is odd,  $N_\a$ contains a Sylow $2$-subgroup of $N$, where $\a\in V$. Of course, $N_\a\ne 1$.
Since $T_i$ is normal in $N$ and $N$ is transitive on $V$,
all $T_i$-orbits on $V$ have equal length, which is a divisor of $|V|$.
Thus $(T_i)_\a$ contains a Sylow $2$-subgroup of $T_i$.

Suppose that $k>1$. Then $T_1$ is intransitive on $V$; otherwise, $T_2$ is semiregular   or acts trivially on $V$, a contradiction.
Assume that $N$ is transitive on $E$. Since $|V|$ is odd, by Lemmas \ref{p-divisor-arc} and \ref{p-divisor-half}, the quotient $\Ga_{T_1}$
is a cycle of odd length. Thus $N/K$ is soluble, where $K$ is the kernel of
$N$ acting on the set of $T_1$-orbits. Then $N=K$, so $N$ fixes every $T_1$-orbits, a contradiction.   Therefore, by Corollaries \ref{odd-order-p-divisor-arc} and \ref{odd-order-p-divisor-half},
$N_\a$ is a $2$-group and $\Ga$ is $G$-arc-transitive.

Let $r=2$. Then $\Ga$ has valency $4$. Again consider the quotient  $\Ga_{T_1}$.
By (\ref{eq-3}),   $\Ga_{T_1}$ has valency no more than $4$.
On the other hand, $\Ga_{T_1}$ is a regular graph of odd order, either
$\Ga_{T_1}$ is a cycle or $\Ga_{T_1}$ has valency $4$. Noting that $(T_1)_\a\ne 1$,  by Lemma \ref{N-semiregular}, $\Ga_{T_1}$ is a cycle of odd length. This leads to a similar contradiction as in the previous paragraph.

Assume that $r$ is odd. Note that $N_\a$ is a Sylow $2$-subgroup of $N$ and
$(T_i)_\a$ is a  Sylow $2$-subgroup of $T_i$; in particular, \[N_\a=(T_1)_\a\times \cdots\times (T_k)_\a.\]
  Since $\Ga$ is $G$-arc-transitive and $N_\a$ is a nontrivial $2$-group, all $N_\a$-orbits on $\Ga(\a)$ have equal length $2$. Take a subgroup $H$ of $G_\a$ of order $r$.
Then $N\cap H=1$, and $H$ acts transitively on the set of  $N_\a$-orbits on $\Ga(\a)$.
Set $X=N{:}H$. Then $\Ga$ is $X$-arc-transitive, and $X_\a=N_\a{:}H$.

Let $M=\l T_1^h\mid h\in H\r$. Then $M$ is a normal subgroup of $X$, $M\le N$ and $M_\a$ is a nontrivial $2$-group. If $M\ne N$ then $M$ is intransitive on $V$ and, by Lemmas \ref{p-divisor-arc}, the quotient $\Ga_M$ is a cycle, which is impossible. Thus $M=N$. Noting that $H$ has prime order $r$,  $H$ acts regularly on $\{T_1,\ldots,T_k\}$ by conjugation; in particular,   $k=r$.
For $h\in H$, setting $T_i^h=T_j$, we have \[(T_i)_\a^h=(T_i\cap N_\a)^h=T_j\cap N_\a=(T_j)_\a.\] It follows that $H$ acts regularly on $\{(T_1)_\a,\ldots,(T_k)_\a\}$ by conjugation.

Let $\b\in \Ga(\a)$. Then $2r=|X_\a:X_{\a\b}|$,
and thus $X_{\a\b}$ is a $2$-group. Since $X_\a$ has a unique Sylow $2$-subgroup,   $X_{\a\b}$ lies in $N_\a$, and so   $X_{\a\b}=N_{\a\b}$ and $|N_\a:N_{\a\b}|=2$.
Since $\Ga$ is connected, by lemma \ref{coset},
there is a $2$-element $x\in \N_X(N_{\a\b})$ with \[\l x,X_\a\r=X,\,
x^2\in N_{\a\b}.\]
Note that $|N\l x\r|={|N||\l x\r:(N\cap \l x \r)|}$ and $|X|=|N||H|$.
We know that $|\l x\r:(N\cap \l x \r)|$ is a divisor of $|H|=r$.
It follows that $|\l x\r:(N\cap \l x \r)|=1$, and so $x\in N$.
Write
$x=x_1x_2\cdots x_k$ with $x_i\in T_i$ for all $i$. Then
\begin{equation}\label{eq-4}
X=\l x,X_\a\r\le \l x_i,(T_i)_\a, H\mid 1\le i\le k \r.
\end{equation}


Consider the projections \[\phi_i: N=T_1\times\cdots\times T_k \rightarrow,\,
t_1t_2\cdots t_k\mapsto t_i.\] Setting $L_i=\phi_i(N_{\a\b})$ for $1\le i\le r$, we have $x_i\in \N_{T_i}(L_i)$, $x_i^2\in L_i$ and  $N_{\a\b}\le L_1\times \cdots L_k$.
Since $N_{\a\b}\le N_\a=(T_1)_\a\times \cdots\times (T_k)_\a$, we have $L_i\le (T_i)_\a$ for $1\le i\le k$.
Recalling that $|N_\a:N_{\a\b}|=2$, we have that $L_i=(T_i)_\a$   for all but at most one of $i$.
Without loss of generality, we let $L_i=(T_i)_\a$   for $i>1$, and  $|(T_1)_\a:L_1|\le 2$.

Since $(T_i)_\a$ is a Sylow $2$-subgroup of $T_i$ and $x_i$ is a $2$-element, if $i>1$ then $x_i\in L_i=(T_i)_\a$ as $(T_i)_\a$ is the unique Sylow $2$-subgroup
of $\N_{T_i}(L_i)$. Recalling that $H$ acts regularly on $\{(T_1)_\a,\ldots,(T_k)_\a\}$ by conjugation, by (\ref{eq-4}), we have
\[X=\l x,X_\a\r\le  \l x_1,(T_1)_\a, H \r.\] Set $H=\{1,h_2,\ldots, h_k\}$. (Note that $k=r$.) Then  \[N{:}H=X=(\l x_1,(T_1)_\a\r\times \l x_1,(T_1)_\a\r^{h_2}\times \cdots\times \l x_1,(T_1)_\a\r^{h_k}){:}H.\]
This implies  that $T_1=\l x_1,(T_1)_\a\r$.

Recall that $|(T_1)_\a:L_1|\le 2$. Then $L_1$ is normal in $(T_1)_\a$.
Noting that $x_1\in \N_{T_1}(L_1)$, it follows that $L_1$ is normal in $\l x_1,(T_1)_\a\r=T_1$, and hence $L_1=1$ as $T_1$ is nonabelian simple and $L_1$ is a $2$-group. Then $|(T_1)_\a|\le 2$, which is impossible as $(T_1)_\a$ is a Sylow  $2$-subgroup of $T_1$. This completes the proof.
\qed

\vskip 5pt

We end this section by a consequence from Theorem \ref{odd-order-prim-2r}.

\begin{corollary}\label{p-power-case}
Assume that $|V|=p^k$ for some odd prime $p$. If $G$ is quasiprimitive on $V$ then
$\Ga$ is arc-transitive, and one of the following holds.
\begin{enumerate}
\item[(1)] $p=3$, $k$ is an odd prime and $\Ga$ is the complete graph $\K_{3^k}$.
\item[(2)] $G\cong \PSU_4(2)$ or $\PSU_4(2).2$, and $\Ga$ has order $27$ and valency $10$.
\item[(3)] $\soc(G)\cong \ZZ_p^k$ and $k\le r$.
\end{enumerate}
\end{corollary}
\proof
Let $N=\soc(G)$. By Theorem \ref{odd-order-prim-2r}, either $N$ is nonabelian simple, or $N\cong \ZZ_p^k$ with $k\le r$.

Assume that $N$ is nonabelian simple.
Note that $|N:N_\a|=p^k$ for $\a\in V$. By \cite{Guralnick}, either $N$ is $2$-transitive on $V$, or $N\cong \PSU_4(2)$ and $N_\a$ has exactly three orbits on $V$ with length $1$, $10$ and $16$ respectively. For these two cases, $\Ga$ is $N$-arc-transitive, and the former case yields (1)  while the latter case gives (2).

Now let $N\cong \ZZ_p^k$. Then we may write $\Ga=\Cay(N,S)$, and $G_\a$ is faithful on $S$, where $\a$ is the vertex corresponding to the identity of $N$. Note that $\{x,x^{-1}\}^h=\{x^h,(x^h)^{-1}\}$ for $x\in S$ and $h\in G_\a$. Then  $G_\a$ induces a transitive action on $\{\{x,x^{-1}\}\mid x\in S\}$. Let $\s:N\rightarrow N,\, y\mapsto y^{-1}$. Then $\s$ is an automorphism of $\Ga$ which fixes $\a$. Thus $\l G_\a, \s\r$ is transitive on $S$, and then
$\Ga$ is $\l G,\s\r$-arc-transitive.
\qed

\vskip 30pt

\section{Proof of Theorem \ref{myth-2}}

Let $\Ga=(V,E)$ be a connected graph of odd order and valency $2r$ for some prime $r$. Let  $G\le \Aut\Ga$ be such that $G$ is primitive on  $V$, transitive on $E$ 
but  not transitive on the $2$-arc set of $\Ga$.

Since $G$ is a primitive group of odd degree, $G$ is known by \cite{Primitive}.
Thus, in this section, we give some further information  about the graph $\Ga$ or its automorphism group. If $\Ga$ is a complete graph, then $G$ is $2$-homogeneous on $V$, and the following result holds.

\begin{lemma}\label{complete-graph}
Assume that $\Ga$ is the complete graph $\K_{2r+1}$ of valency $2r$. Then
\begin{enumerate}
\item[(1)] $G\le \AGammaL_1(p^d)$, $r={p^d-1\over 2}$, and either $d=1$ or $p=3$ and $d$ is an odd prime; or
    \item[(2)] $\ASL_d(3)\le G\le \AGL_d(3)$, and $d$ is an odd prime; or
      \item[(3)] $G=\PSL_2(11)$, $\A_7$ or $\PSL_d(2)$, and
      $r=5$, $7$ or $2^{d-1}-1$ respectively.
\end{enumerate}
\end{lemma}
\proof
Note that $G$ is $2$-homogeneous on $V$. Since $\Ga$ is not $(G,2)$-arc-transitive, $G$ is not $3$-transitive on $V$.

Suppose that $G$ is an affine primitive group.
Let $G\le \AGL_d(p)$. Since $p^d=|V|=2r+1$, either $d=1$ or
$p=3$ and $d$ is an odd prime. If $G\not\lesssim \AGammaL_1(p^d)$ then, by \cite{Kantor}, $G$ is $2$-transitive, and then $G\ge \ASL_d(p)$ by checking
\cite[Table 7.3]{Cameron}. Thus
(1) or (2) occurs.

Suppose that $G$ is almost simple. Then $G$ is $2$-transitive but not  $3$-transitive. Checking
\cite[Table 7.4]{Cameron}, we get (3).
\qed

\begin{lemma}\label{over-group-2-arc}
 Let $G\le X\le \Aut\Ga$  and $\a\in V$. Suppose that $\Ga$ is $(X,2)$-arc-transitive.
 Then one of the following holds.
 \begin{enumerate}
 \item[(1)] $X_\a^{\Ga(\a)}=\A_{4}$ or $\S_4$, and $r=2$;
 \item[(2)] $\soc(X_\a^{\Ga(\a)})=\A_{2r}$, $r\ge 3$;
 \item[(3)] $\soc(X_\a^{\Ga(\a)})=\PSL(2,q)$, and $q=5$ or $q=2r-1\ge 9$;
 \item[(4)]  $\soc(X_\a^{\Ga(\a)})=\M_{22}$, and $r=11$.
 \end{enumerate}
 \end{lemma}
 \proof
Since $X_\a^{\Ga(\a)}$ is a $2$-transitive group of degree $2r$,
the lemma follows from checking the finite $2$-transitive groups,
 refer to \cite[Tables 7.3 and 7.4]{Cameron}.
 \qed


 \begin{lemma}\label{over-group-1-arc}
 Let $G\le X\le \Aut\Ga$  and $\a\in V$. Suppose that $X_\a$ is insoluble and $\Ga$ is not $(X,2)$-arc-transitive.
 Then $X_\a$ has a composition factor isomorphic to one of the following simple groups.
 \begin{enumerate}
\item[(1)] $\PSL(2,11)$, and $r=11$;
\item[(2)] $\M_{11}$ or $\M_{23}$, and $r=11$ or $23$, respectively;
\item[(3)] $\PSL(d,q)$, and $r={q^d-1\over q-1}$, where $q$ is prime power and $d$ is  a prime;
    \item[(4)] $\A_r$ for $r\ge 5$.
\end{enumerate}
 \end{lemma}
 \proof
 Let $\Delta(\a)$ be an $X_\a$-orbit on $\Ga(\a)$. Then $|\Delta(\a)|\in \{r,2r\}$, and $X_\a^{\Delta(\a)}$ is insoluble, refer to \cite[Theorem 3.2C]{Dixon}. If $|\Delta(\a)|=r$ then $X_\a^{\Delta(\a)}$ is a $2$-transitive group of prime degree, and the lemma follows from checking the finite $2$-transitive groups.

 Assume that $\Delta(\a)=\Ga(\a)$. If $X_\a^{\Ga(\a)}$ is primitive then, by \cite[Theorem 1.51]{Gorenstein}, $r=5$ and $X_\a^{\Ga(\a)}=\A_5$ or $\S_5$,
 and thus the lemma follows.
 Suppose next that  $X_\a^{\Ga(\a)}$ is imprimitive. Then $X_\a^{\Ga(\a)}\le \S_2\wr\S_r$ or $\S_r\wr\S_2$.

 Let $X_\a^{\Ga(\a)}\le \S_2\wr\S_r$. Let $O$ be the largest normal $2$-subgroup of the wreath product $\S_2\wr\S_r$. Then $\S_r\cong \S_2\wr\S_r/O\ge X_\a^{\Ga(\a)}O/O\cong X_\a^{\Ga(\a)}/(X_\a^{\Ga(\a)}\cap O)$. Since $X_\a^{\Ga(\a)}$ is insoluble,  $X_\a^{\Ga(\a)}/(X_\a^{\Ga(\a)}\cap O)$ is isomorphic to a $2$-transitive subgroup of $\S_r$, and then the lemma follows.

 Let $X_\a^{\Ga(\a)}\le  \S_r\wr\S_2$. Then $X_\a^{\Ga(\a)}\cap \S_r^2$ has
 two orbits on $\Ga(\a)$. Considering the action of $X_\a^{\Ga(\a)}\cap \S_r^2$ on one of its orbits, we get the lemma.
 \qed

\begin{theorem}\label{over-group}
Let $G\le X\le \Aut\Ga$. Then either $\soc(G)=\soc(X)$, or $\Ga \cong \K_{2r+1}$.
\end{theorem}
\proof
Assume that $\Ga$ is not a complete graph. Then $X$ is not $2$-homogeneous on $V$.
We next suppose that $\soc(G)\ne \soc(X)$, and deduce the contradiction.

Suppose that $\soc(G)\ne \soc(X)$.
By Theorem \ref{myth} and \cite{Li-odd}, $X$ is either affine or almost simple if
$X$ is  not transitive on the $2$-arc set of $\Ga$, and  $X$  is almost simple if $\Ga$ is $(X,2)$-arc-transitive.

Assume that either  $G$ or $X$ is an affine primitive group.
By \cite[Propositions 5.1 and 5.2]{Praeger-inclusion}, we have $G\le \AGL_3(3)$ and $\soc(X)=\PSU_4(2)$.
Then $X_\a\cong \ZZ_2^4{:}\A_5$ or $\ZZ_2^4{:}\S_5$ for $\a\in V$, and so
$r=5$ by Lemma \ref{normal-subg-1}; however, $5$ is not a divisor of $|G|$, a contradiction.

Assume next that both $G$ and $X$ are almost simple. Let $\a\in V$. Recall  that $|V|$ is odd and $X$ is not $2$-homogeneous on $V$. By
\cite[Proposition 6.1]{Praeger-inclusion} and \cite{LPS}, all possible triples
$(\soc(G),\soc(X), |V|)$ are listed in  Table \ref{over-almost}.
\begin{table}
\[\begin{array}{|c|c|c|c|c|} \hline
\mbox{Line}&\soc(G) &\soc(X) & |V| & \mbox{$X$-action, remark} \\  \hline
1&\M_{11}&\A_{11}& 55,165&   \\ \hline
2&\M_{12}&\A_{12}& 495&   \\ \hline
3&\M_{22}&\A_{22}& 231&   \\ \hline
4&\M_{23}&\A_{23}& 253,1771&   \\ \hline
5&\PSL_2(q)&\A_{q+1}& {1\over 2}q(q+1)& 2\mbox{-sets},  q  \equiv 1\, (\mod 4)\\ \hline
6&\A_{2l-1}&\A_{2l}& {1\over 2}\binom{2l}{l}& l,l-\mbox{partitions}, l\ge 3\\ \hline
7&\PSL_2(11)&\M_{11} &55  &   \\ \hline
8&\M_{23}&\M_{24} &1771  &   \\ \hline
9&\PSL_m(3)&\POmega^+_{2m}(3)&{1\over 2}3^{m-1}(3^m-1)&\mbox{nonsingular points}, m \mbox{ odd} \\  \hline
10&\G_2(q)&\Omega_7(q)&{1\over 2}q^3(q^3\pm 1)& \mbox{nonsingular
hyperplanes}\\
&&&&q  \equiv \pm 1\, (\mod 4)\\  \hline
11&\A_{12}& \POmega^-_{10}(2)& 495&  \\  \hline
12&\J_3&\PSU_9(2)& 43605&   \\  \hline
13&\Omega_7(3)& \POmega^+_8(3)& 28431&  \\  \hline
14&\G_2(3)& \Omega_7(3)&3159&  \\  \hline
\end{array}\]
{\caption{}\label{over-almost}}
\end{table}

For Line 10 of  Table \ref{over-almost}, by \cite[Table 8.39]{Low-dim}, $X_\a$ has a unique insoluble composition factor, say $\POmega^\pm_6(q)$, which contradicts Lemmas \ref{over-group-2-arc} and \ref{over-group-1-arc}.
For Line 13 or 14 of  Table \ref{over-almost}, by the Atlas \cite{Atlas}, $X_\a$ is an almost simple group with socle $\POmega^+_8(2)$ or $\Sp_6(2)$ respectively, and we get a similar contradiction.

Note that $r\in \pi(G_\a)\cap \pi(X_\a)$ and, by Lemma \ref{normal-subg-1},
$\max(X_\a)<2r$. This allows us exclude Lines 1-4, 11, and 12 of Table \ref{over-almost}. For example, if Line 12 of  Table \ref{over-almost} occurs then $G_\a$ is a $\{2,3\}$-group by the Atlas \cite{Atlas}, so $r\le 3$, yielding that
$|X_\a|$ has no prime divisor other than $2,\,3$ and $5$, which is impossible as
$|X_\a|$ is divisible by $11$.

Note that $\Ga(\a)$ is either an $X_\a$-orbit of length $2r$ or
the union of two  $X_\a$-orbits of length $r$. In particular, $X_\a$ has a subgroup of index a prime or twice a prime.
For Lines 7 and 8 of  Table \ref{over-almost}, the lengthes of $X_\a$-orbits on $V$ are known, refer to the Webpage Edition of \cite{Atlas}.
If Line 7 of  Table \ref{over-almost} occurs, then $X$ is a primitive group (on $V$) of rank $3$, and $X_\a$ has three orbits on $V$ with length $1$, $18$ and $36$ respectively, which gives a contradiction. For Line 8 of  Table \ref{over-almost}, $X_\a$ has four orbits on $V$ with length $1$, $90$, $240$ and $1440$ respectively, we get a similar contradiction.
For Line 9 of  Table \ref{over-almost}, by \cite{HN}, $X_\a$ has three orbits on $V$ with length $1$, ${1\over 2}3^{m-1}(3^{m-1}-1)$ and $3^{2m-2}-1$ respectively; however, none of these three numbers has the form of $2r$ or $r$.

Suppose that $\soc(G)$ and $\soc(X)$ are given as in Line 6 of Table \ref{over-almost}. Then action of $G$ on $V$ is equivalent that on the set
of $(l-1)$-subsets of $\{1,2,\ldots ,2l-1\}$. Thus $G_\a$ has $l$ orbits on $V$,
which have lengthes $\binom{l-1}{i}\binom{l}{l-1-i}$, $0\le i\le l-1$, respectively.
If $r=\binom{l-1}{i}\binom{l}{l-1-i}$ for some $i$, then $i=0$ and $l=r$.
If $2r=\binom{l-1}{i}\binom{l}{l-1-i}$ for some $i$, then either $i=0$ and $l=2r$ or $i=1$ and $l=3$. For each of these three cases, it is easily shown that
${1\over 2}\binom{2l}{l}$ is even, which is not the case as $|V|$ is odd.

Finally, let $\soc(G)$ and $\soc(X)$ be as in Line 5 of Table \ref{over-almost}.
Then $X_\a$ has three orbits on $V$ with length $1$, $2(q-1)$ and ${1\over 2}(q-1)(q-2)$, respectively. The only possibility is that $q=5$,  $\soc(G)=\PSL_2(5)$ and $2r=6$;  however,
in this case, $G_\a$ is a $2$-group which is not maximal in $G$, a contradiction.
This completes the proof.
\qed

\begin{remark}
For Line 8 of  Table \ref{over-almost}, one may construct a graph of valency $90$ and order $1771$ which is both $M_{23}$-arc-transitive and $M_{24}$-arc-transitive. Thus Theorem \ref{over-group} does not hold without the assumption that $\Ga$ has twice prime valency.
\end{remark}

By Theorem \ref{over-group}, we have  the following consequence which finishes the proof of Theorem \ref{myth-2}.

\begin{corollary}\label{cor-1}
$\Ga$ is $2$-arc-transitive if and only if $\Ga$ is a complete graph.
\end{corollary}
\proof
Note that a complete graph must be $2$-arc-transitive. Thus, it suffices to show that $\Ga$ is not $2$-arc-transitive if $\Ga\not\cong\K_{2r+1}$.

Suppose that $\Ga\not\cong\K_{2r+1}$ and  $\Ga$ is  $2$-arc-transitive. Let $N=\soc(G)$ and $X=\Aut\Ga$. Then $N=\soc(X)$ by Theorem \ref{over-group}.
Noting that $|V|$ is odd, by Theorem \ref{odd-order-prim-2r} and \cite{Ivanov-Praeger-93}, $N$ is a nonabelian simple group.

Let $\a\in V$. Since $|V|$ is odd, $N_\a\ne 1$, and so $N_\a^{\Ga(\a)}\ne 1$.
Recalling that $\Ga$ is not $(G,2)$-arc-transitive, $G_\a^{\Ga(\a)}$ and hence $N_\a^{\Ga(\a)}$ is not $2$-transitive. Noting that $N_\a^{\Ga(\a)}$
is a proper normal subgroup of $X_\a^{\Ga(\a)}$,   by Lemma \ref{over-group-2-arc}, we have $r=2$, $N_\a^{\Ga(\a)}\cong \ZZ_2^2$.
In particular, $\Ga$ is a primitive $2$-arc-transitive graph of valency $4$.
It follows from  \cite[Theorem 1.5]{LLM} that $|N_\a|$ is divisible by $3$.
Then, by Lemma \ref{normal-subg-1},  $N_\a^{\Ga(\a)}$ is not a $2$-group, a contradiction.
This completes the proof.
\qed

We end this section by another consequence of Theorem \ref{over-group}.

\begin{corollary}\label{cor-2}
Assume that $\Ga \not\cong \K_{2r+1}$. Then
either $\Ga$ is $G$-arc-transitive, or
$r\ge 3$ and $\Aut\Ga$ has a subgroup of index at most $2$ which is transitive on the edge set but not transitive  on the arc set of $\Ga$.
\end{corollary}
\proof
Since
$G$ is primitive on $V$ and $\Ga$ is not a cycle, for $\a\in V$,
each $G_\a$-orbit on $V$ has length at least $3$.
Thus, if $r=2$ then $\Ga$ is $G$-arc-transitive.

Suppose that $\Ga$ is not $G$-arc-transitive.
Then $r\ge 3$, and $G$ has two orbits on the arc set of $\Ga$, say $\Delta$ and $\Delta^*=\{(\b,\a)\mid (\a,\b)\in \Delta\}$. To complete the proof of this corollary, we assume further that $\Ga$ is arc-transitive. Let $X=\Aut\Ga$.
By Theorem \ref{over-group}, $N:=\soc(G)=\soc(X)$.
Noting that $N$ fixes both $\Delta$ and $\Delta^*$ set-wise,
$X$ induces a transitive action on $\{\Delta,\Delta^*\}$. Then the result follows.
\qed

\vskip 30pt

\section{Examples and a proof of   Theorem \ref{myth-3}}
We first list several examples of graphs which may support some results in Section 4.

\vskip 5pt

Let $G$ be a finite group and $H$ be a subgroup of $G$ with $\cap_{g\in  X}H^g=1$.
Then $G$ acts (faithfully) on the set $[G:H]$ of right cosets of $H$ in $G$ by \[y: Hx\mapsto Hxy, \forall x,y\in G.\] Take a $2$-element $x\in G\setminus H$ with $x^2\in H$.
Then $x$ normalizes $K:=H\cap H^x$. Define a graph $\Cos(G,H,x)=(V,E)$ with
\[V=[G:H],\, E=\{\{Hg_1,Hg_2\}\mid g_1,g_2\in G,\,g_2g_1^{-1}\in HxH\}.\]
Then $\Cos(G,H,x)$ is a $G$-arc-transitive graph of order $|G:H|$ and valency $|H:K|$. Moreover, it is well-known that the graph $\Cos(G,H,x)$ is connected if and only if $\l x, H\r=G$.
The first example gives three graphs which meet Corollary \ref{odd-order-p-divisor-arc}.
\begin{example}\label{exam-1}
{\rm
(1) Let $T=\PSL(2,27)$ and $G=T{:}\l \s\r$, where $\s$ is induced by the Frobenius  automorphism of the field of order $27$. Then, by the Atlas \cite{Atlas}, $G$ has a subgroup $H$ with $H\cong \A_4$ and $H\cap T\cong \ZZ_2^2$.

Let $o\in H\cap T$. Then $\C_T(o)\cong \D_{28}$, and $\C_T(o)\cap H=H\cap T\cong \ZZ_2^2$. Confirmed by the GAP \cite{GAP}, there is an involution $x\in \C_T(o)\setminus H$ such that $|H\cap H^x|=2$ and $\l x,H\r=G$.  Then $\Ga=\Cos(G,H,x)$  is a $G$-arc-transitive graph of valency $6$ and order $2457$.

Clearly, $T$ is transitive on the vertex set but not transitive on the edge set of $\Ga$. In fact, $\Ga$ is the edge-disjoint union of three
$T$-arc-transitive graphs of valency $2$, and each of them is the vertex-disjoint union of $351$ copies of the $7$-cycle.

\vskip 5pt

(2) Let $G=\PGL(2,11)$, $T=\PSL(2,11)$, $H\cong \D_{24}$ and $H_1=T\cap H$.
Let $K\le H_1$ with $K\cong \ZZ_2^2$. Then $\N_G(K)\cong \S_4$, $\N_H(K)\cong \D_8$ and $\N_T(K)\cong \A_4$.

Take an involution $x\in \N_G(K)\setminus H$.
Then $\l x, H\r=X$ and $H\cap H^x=K$.
Thus $\Cos(G,H,x)$ is a connected $G$-arc-transitive graph of valency $6$ and order $55$. It is easily to see that $T$ acts transitively on the edges but not on the arcs of  $\Cos(G,H,x)$.

\vskip 5pt

(3) Let $G=\PGL(2,17)$, $T=\PSL(2,17)$ and $H=\l c\r \times H_1$, where $T\ge H_1\cong \S_3$ and
$c\in G\setminus T$ is an involution. Let  $o \in H_1$ be an involution.
Then $\C_T(o)\cong \D_{16}$.

Take an involution $x\in \C_T(o)\setminus H$. Then
$\l x, H\r=G$, $\l x, H_1\r\cong \S_4$, $|H:(H\cap H^x)|=6$ and  $|H_1:(H_1\cap H_1^x)|=6$, confirmed by the GAP \cite{GAP}. It follows that
$\Cos(G,H,x)$ is a connected $G$-arc-transitive graph of valency $6$ and order $408$.
Moreover, $\Cos(G,H,x)$ is the edge-disjoint union of two
$T$-arc-transitive graphs of valency $3$, and each of them is the vertex-disjoint union of $102$ copies of the complete graph $\K_4$.
}\qed
\end{example}

\begin{example}\label{exam-2}
{\rm (1) Let $G=\PGL(2,7)$, $T=\PSL(2,7)$ and $\D_{16}\cong H\le G$.
Let $\ZZ_2^2\cong K\le H\cap T$. Then $N_H(K)\cong \D_8$ and $N_G(K)=N_T(K)\cong \S_4$.

Take an involution $x\in N_T(K)\setminus H$. Then $\Ga=\Cos(G,H,x)$ is a connected
$G$-arc-transitive graph of valency $4$ and order $21$. The graph $\Ga$ is in fact the edge-disjoint union of two
$T$-arc-transitive graphs of valency $2$, and each of them is the vertex-disjoint union of $7$ copies of the $3$-cycle.

\vskip 5pt

(2) Let $G=\PGL(2,9)$, $T=\PSL(2,9)$ and $\D_{16}\cong H\le G$.
Take $\ZZ_2^2\cong K\le H\cap T$. Then $N_H(K)\cong \D_8$ and $N_G(K)=N_T(K)\cong \S_4$. For  an involution $x\in N_T(K)\setminus H$, we have a $G$-arc-transitive graph  $\Ga=\Cos(G,H,x)$
of valency $4$ and order $45$. Further, the graph $\Ga$ is the edge-disjoint union of two
$T$-arc-transitive graphs of valency $2$, and each of them is the vertex-disjoint union of $15$ copies of the $3$-cycle.
}\qed
\end{example}

\vskip 5pt

Note that the graphs in Example \ref{exam-2} are all vertex-primitive. By \cite{LLM}, up to isomorphism, there are exactly two vertex-primitive arc-transitive graphs of valency $4$ which are not $2$-arc-transitive, and then the following result holds.

\begin{lemma}\label{prim-4}
Let $\Ga=(V,E)$ be a connected graph of valency $4$, and let $G\le \Aut\Ga$.
Suppose that $G$ is almost simple and transitive on $E$.
If $G$ is primitive on $V$ then either $\Ga$ is $\soc(G)$-arc-transitive, or one of the following holds:
\begin{enumerate}
\item[(1)] $\Aut\Ga=G=\PGL(2,7)$ and $\Ga$ is given in (1) of Example \ref{exam-2};
\item[(2)] $G=\PGL(2,9)$, $\M_{10}$ or $\PGammaL(2,9)$, and $\Ga$ is given in (2) of Example \ref{exam-2}.
\end{enumerate}
\end{lemma}

\vskip 5pt

\noindent {\bf Proof Theorem \ref{myth-3}}.
Let $\Ga=(V,E)$ be a connected    graph of odd order and twice prime valency $2r$, and let $G\le \Aut\Ga$ be such that $\Ga$ is $G$-edge-transitive but not $(G,2)$-arc-transitive. Assume that $G$ is almost simple and primitive on $V$.

Let $N=\soc(G)$. Note that $N$ is transitive but not regular on $V$.
If $\Ga$ is not $G$-arc-transitive then, since $G$ is primitive,
$r$ is odd in this case, and $\Ga$ is $N$-edge-transitive by Corollary \ref{odd-order-p-divisor-half}. Thus we assume that $\Ga$ is  $G$-arc-transitive.
Then, by
Corollary \ref{odd-order-p-divisor-arc}, either $\Ga$ is $N$-edge-transitive
or $N_\a$ is a Sylow $2$-subgroup of $N$, where $\a\in V$.

Assume that $N_\a$ is a Sylow $2$-subgroup of $N$.
Since $G$ is a primitive group of odd degree, by \cite{Primitive},
we conclude that $N\cong \PSL(2,q)$ and $N_\a$ is a dihedral group.
It follows from  \cite[Table 8.1]{Low-dim} that
$N_\a\cong \D_{q-1}$ or $\D_{q+1}$.

Suppose that $|N_\a|>8$. Consider the Frattini subgroup $\Phi$ of $N_\a$.
Clearly, $\Phi\ne 1$ is a characteristic subgroup of $N_\a$, and $\Phi$ is contained in the unique cyclic subgroup of $N_\a$ with index $2$. It follows that
$\Phi$ is a characteristic subgroup of every maximal subgroup of $N_\a$.
Take $R\le G_\a$ with $|R|=r$, set $X=N{:}R$. Then $\Ga$ is $X$-arc-transitive, $X_\a=N_\a{:}R$ and $X_{\a\b}=N_{\a\b}$, where $\b\in \Ga(\a)$.
Since $\Ga$ is connected, by Lemma \ref{coset}, $X=\l x,N_\a R\r$ for some $2$-element $x\in \N_X(N_{\a\b})$. It is easily shown that $x\in N$, and so
$x\in N_N(N_{\a\b})$.
Noting that $|N_\a:N_{\a\b}|=2$, it follows that $\Phi$ is  a characteristic subgroup of $N_{\a\b}$. Then both $R$ and $x$ normalize $\Phi$, and so
$\Phi$ is normal in $X=\l x,N_\a R\r$, which is impossible.

Finally, let $|N_\a|\le 8$. Then we have $q=7$ or $9$, and $N\cong \PSL(2,7)$ or $\PSL(2,9)$, respectively. Checking the maximal subgroups of $G$, we conclude that $G_\a$ is a (Sylow) $2$-subgroup of $G$. Thus we have $r=2$, and $\Ga$ has valency $4$. Then our theorem follows from Lemma \ref{prim-4}.


\vskip 50pt

\begin{thebibliography}{99} 



\bibitem{Low-dim} J.N. Bray, D.F. Holt and C.M. Roney-Dougal,
\emph{The maximal subgroups of the low-dimensional finite classical
groups}, Cambridge University Press,  2013.

\bibitem{Cameron} P.J. Cameron, {\em  Permutation Groups}, Cambridge University Press,  Cambridge, 1999.


\bibitem{CLP2000} M.D. Conder, C.H. Li and C.E. Praeger,
On the weiss conjecture for finite locally primitive graphs, \emph{
P. Edinburgh Math. Soc.}  {\bf 43} (2000), 129-138.

\bibitem{Atlas}
J.H. Conway, R.T. Curtis, S.P. Noton, R.A. Parker and R.A. Wilson,
{\em Atlas of Finite Groups}, Clarendon Press, Oxford, 1985. \\
{\tt{http://brauer.maths.qmul.ac.uk/Atlas/v3/}}

\bibitem{Dixon}
 J.D. Dixon and B. Mortimer, {\em Permutation Groups},
Springer-Verlag  New  York  Berlin  Heidelberg,  1996.


\bibitem{Guralnick}
R. M. Guralnick, Subgroups of prime power index in a simple group, {\em J. Algebra} {\bf 81} (1983), 304-311.






\bibitem{Gorenstein} D. Gorenstein, {\em Finite Simple Groups},
Pleunum Press, New Youk, 1982.

\bibitem{HN} J.I. Hall and H.N. Nguyen,
On thepermutation modules for orthogonal groups ${\rm O}^\pm_m(3)$ acting on nonsingular points of their standardmodules, {\em J. Pure Appl. Algebra}
{\bf 216} (2012), 581-592.




\bibitem{Ivanov-Praeger-93}
A.A. Ivanov and C.E. Praeger,
On finite affine $2$-arc transitive graphs, {\em Europ. J. Combin.}
{\bf 14} (1993), 421-444.


\bibitem{Kantor}W.M. Kantor, $k$-homogeneous groups, {\em Math. Z.} {\bf 124} (1972), 261-265.




\bibitem{Li-odd} C.H. Li, On finite s-transitive graphs of odd order,
 {\em J. Combin. Theory Ser. B} {\bf 81} (2001), 307-317.

\bibitem{LPS} M.W. Liebeck, C.E. Praeger and J. Saxl,
A classification of the maximal subgroups of the finite alternating
and symmetric groups, {\em J. Algebra }  {\bf 111}(2) (1987),
365-383.



\bibitem{LLM} C.H. Li, Z.P. Lu and D. Maru\v{s}i\v{c},
On primitive permutation groups with small suborbits and their
orbital graphs, \emph{J. Algebra} {\bf 279} (2004),749-770.

\bibitem{LLiuL}
C.H.Li, Z. Liu and Z.P. Lu, The edge-transitive tetravalent Cayley graphs of square-free order, {\em Discrete Math.} {\bf 312} (2012), 1952-1967.

\bibitem{Primitive} M.W. Liebeck and J. Saxl, The primitive permutation
groups of odd degree, \emph{J. London Math. Soc.} (2) {\bf 31} (1985),
250-264.

\bibitem{Neumann}
P.M. Neumann, Finite permutation groups, edge-colored graphs and
matrices, {\em Topics in Group Theory and Computation}, Acad. Press,   1977, pp.
82-118.



\bibitem{Praeger-inclusion}
C.E. Praeger, The inclusion problem for finite primitive permutation groups, {\em Proc. London Math. Soc.} (3) {\bf  60} (1990), 68-88.

\bibitem{Prag-o'Nan} C.E. Praeger, An O'Nan-Scott theorem for finite
quasiprimitive permutation groups and an application to $2$-arc
transitive graphs, {\em J. London Math. Soc.} {\bf 47}(1992),
227-239.



\bibitem{Prag-quasi}
C.E. Praeger, Finite quasiprimitive graphs,
 {\em London Math. Soc. Lecture Note Ser.}, vol. {\bf 241}, Cambridge
University Press, 1997, pp. 65-85.


\bibitem{Prag-Xu} C.E. Praeger and M.Y. Xu, A characterization of a class of symmetric graphs of twice prime valency, {\em Europ. J. Combin.} {\bf 10} (1989), 91-102.


\bibitem{Sabiddusi} G.O. Sabiddusi,
Vertex-transitive graphs, {\em  Monatsh. Math.}, {\bf 68}(1964),
426-438.



\bibitem{GAP} The GAP Group, GAP-Groups, Algorithms, and Programming, Version 4.8.6, 2017. http://www.gap-system.org





\end{thebibliography}
\end{document}